\documentstyle[12pt,leqno]{article}
\begin{document}
\title{Nambu-Lie Groups}
\author{{\normalsize by}\\Izu Vaisman}
\date{}
\maketitle
{\def\thefootnote{*}\footnotetext[1]%
{{\it 1991 Mathematics Subject Classification}
58 F 05, 22 E 99. \newline\indent{\it Key words and phrases}:
Nambu brackets, Nambu-Lie groups, Nambu-Lie algebras.}}
\begin{center} \begin{minipage}{12cm}
A{\footnotesize BSTRACT.
We extend the Nambu bracket
to $1$-forms.
Following the Poisson-Lie case, we define
Nambu-Lie groups as Lie groups endowed with a multiplicative
Nambu structure.
A Lie group $G$ with a Nambu structure $P$
is a Nambu-Lie group iff $P=0$ at the unit, and the
Nambu bracket of left (right) invariant forms
is left (right) invariant. We define a corresponding notion of a
Nambu-Lie algebra.
We give several examples of Nambu-Lie groups and algebras.}
\vspace{5mm}
\end{minipage} \end{center}
\indent
In 1973, Nambu \cite{Nb} studied a dynamical system which was defined as a
Hamiltonian system with respect to a ternary,
Poisson-like bracket defined by a
Jacobian determinant.
A few years ago, Takhtajan \cite{Tk1} reconsidered the
subject, proposed a general, algebraic definition of a {\it
Nambu-Poisson bracket of order} $n$ which, for brevity,
we call a {\it Nambu
bracket},
and gave the basic
properties of this operation. The Nambu bracket
is an intriguing operation, in spite of its rather restrictive
character, which follows from the fact
conjectured in \cite{Tk1}, and proven by several authors
\cite{Gt}, \cite{AG},
\cite{Nks}, \cite{Le}, \cite{MVV} namely, that,
locally and with respect to well chosen coordinates,
any nonzero Nambu bracket is just a Jacobian
determinant.

In this paper, we show that a Nambu bracket induces a corresponding
bracket of $1$-forms, and use the latter for a characterization of {\it
Nambu-Lie groups}, a natural generalization of the Poisson-Lie groups
(e.g, \cite{V1}). The relation with a
corresponding notion of a Nambu-Lie algebra is discussed,
and several examples of Nambu-Lie groups and algebras are given.

Independently, the Nambu-Lie groups have been studied by J. Grabowski and
G. Marmo,
who determined the general structure of the multiplicative Nambu
tensor fields on Lie groups \cite{GM}.
A preliminary version of our paper circulated as a preprint before 
\cite{GM} was available. But, in the present, final version we will also use results
from 
\cite{GM}, with due quotation.

{\it Acknowledgements}. The final version of this 
paper was written during visits of
the author at the Department of Mathematics, Pennsylvania State University, U.S.A.,
and the Erwin Schr\"odinger International Institute for Mathematical Physics, Vienna,
Austria. I want to express here my gratitude to these host institutions for their
support, and to  Jean-Luc Brylinski and  Ping Xu at Penn State, and Peter Michor at
ESI-Vienna for their invitations there.
\section{Nambu brackets of
functions and 1-forms} Since the subject is not classical,
we first recall the notion of a Nambu bracket and the
geometric structure behind it.
Let $M^{m}$ be an $m$-dimensional differentiable
manifold (in this paper everything is of the $C^{\infty}$ class), and ${\cal
F}(M)$ its algebra of real valued functions. A {\it
Nambu bracket or structure} of order $n$, $3\leq n\leq m$, is an internal
$n$-ary operation on ${\cal F}(M)$, denoted by $\{\;\}$, which satisfies the
following axioms:\vspace{2mm}\\
$(i)$ $\{\;\}$ is {\bf R}-multilinear and totally
skew-symmetric;\vspace{2mm}\\
$(ii) \hspace{1cm}\{f_{1},\ldots,f_{n-1},gh\}=\{f_{1},\ldots,f_{n-1},g\}h
+g\{f_{1},\ldots,f_{n-1},h\}$\vspace{2mm}\\
(the {\it Leibniz rule});\vspace{2mm}\\
$(iii). \hspace{2cm}\{f_{1},\ldots,f_{n-1},\{g_{1},\ldots,g_{n}\}\}=$
$$\sum_{k=1}^{n}\{g_{1},\ldots,g_{k-1},
\{f_{1},\ldots,f_{n-1},g_{k}\},g_{k+1},\ldots,g_{n}\}$$
(the {\it fundamental identity}).
A manifold endowed with a Nambu bracket is a {\it Nambu
manifold}. (Remember that if we use the same definition for $n=2$, we get a
{\it Poisson bracket}.)

By $(ii)$, $\{\;\}$ acts on each factor as a vector field, whence it must be
of the form
$$\{f_{1},\ldots,f_{n}\}=P(df_{1},\ldots,df_{n}),\leqno{(1.1)}$$
where $P$ is a field of $n$-vectors on $M$. If such a field defines a
Nambu bracket, it is called a {\it Nambu tensor (field)}.
$P$ defines a bundle mapping
$$\sharp_{P}:\underbrace
{T^{*}M\times\ldots\times T^{*}M}_{n-1\:times}\longrightarrow TM \leqno{(1.2)}$$
given by
$$<\beta,\sharp_{P}(\alpha_{1},\ldots,\alpha_{n-1})>=
P(\alpha_{1},\ldots,\alpha_{n-1},\beta) \leqno{(1.3)}$$
where all the arguments are covectors.

In what follows, we denote an $n$-sequence of functions or forms, say
$f_{1},\ldots, f_{n}$, by $f_{(n)}$, and, if an index $k$ is missing, by
$f_{(n,\hat k)}$.

The next basic notion is that of the $P$-{\it Hamiltonian vector
field} of $(n-1)$ functions defined by
$$X_{f_{(n-1)}}=\sharp_{P}(df_{(n-1)}).\leqno{(1.4)}$$
The fundamental identity $(iii)$ means that {\it the Hamiltonian vector
fields are derivations of the Nambu bracket}. Another
interpretation of $(iii)$ is
$$(L_{X_{f_{(n-1)}}}P)(dg_{1},\ldots,dg_{n})=0,\leqno{(1.5)}$$
where $L$ is the Lie derivative, i.e., {\it the Hamiltonian vector fields
are infinitesimal automorphisms of the Nambu tensor}.

A mapping $\varphi:(M_{1},P_{1})
\rightarrow(M_{2},P_{2})$
between two Nambu manifolds of the same order $n$
is a {\it Nambu morphism} if the tensor
fields $P_{1}$ and $P_{2}$
are
$\varphi$-related or, equivalently,
$\forall g_{(n)}\in{\cal F}(M_{2})$, one has
$$\{g_{1}\circ\varphi,\ldots,g_{n}\circ\varphi\}_{1}=
\{g_{1},\ldots,g_{n}\}_{2}.$$
If, moreover, $\varphi$ is a diffeomorphism ,
it is an {\it equivalence of
Nambu manifolds}. The notion of a Nambu morphism also
leads to the following definition: a submanifold $N$ of a Nambu manifold 
$(M,P)$ is a {\it Nambu submanifold} if $N$ has a (necessarily unique)
Nambu tensor $Q$ such that the immersion of $(N,Q)$ in $(M,P)$ is a
Nambu morphism. Like in the Poisson case, $Q$ exists iff, along $N$, $P$ vanishes
whenever evaluated on $n$ $1$-forms one of which,
at least, belongs to the annihilator space 
$Ann(TN)$, and then $im\sharp_{P}$ is a tangent 
distribution of $N$ (e.g., \cite{V1}).
\proclaim 1.1 Theorem. \cite{Gt}, \cite{AG}, \cite{Nks}, \cite{Le},
\cite{MVV}.
$P$ is a Nambu tensor field of order $n$ iff $\forall p\in M$
where $P_{p}\neq0$ there are local coordinates
$(x^{k},y^{\alpha})$ ($k=1,\ldots,n$, $\alpha=1,\ldots,m-n$) around $p$
such that $$P=\frac{\partial}{\partial x^{1}}\wedge\ldots\wedge
\frac{\partial}{\partial x^{n}}\hspace{1cm}
(\{f_{1},\ldots,f_{n}\}=\frac{\partial
(f_{1},\ldots,f_{n})}{\partial
(x^{1},\ldots,x^{n})}).\leqno{(1.6)}$$
\par

On the canonical coordinate neighborhood where
(1.6) holds we have $${\cal
D}:=span(im\,\sharp_{P})=span\{\partial/\partial x^{k}\}.$$
Hence, globally ${\cal D}$
is a foliation with singularities whose leaves are either points,
called the {\it singular points} of $P$,
or $n$-dimensional submanifolds with a Nambu
bracket induced by $P$. ${\cal D}$ is the
{\it canonical foliation} of the Nambu structure $P$. The canonical
foliation is regular i.e., all the leaves are $n$-dimensional, iff
$P$ never vanishes, and then we say that $P$ is a {\it regular
Nambu structure}.

\proclaim 1.2 Theorem. \cite{Gt}, \cite{Le}. A regular Nambu structure of
order $n$ on a
differentiable manifold $M^{m}$ is equivalent with a regular
$n$-dimensio\-nal
foliation $S$ of $M$, and a bracket operation defined by the formula
$$d_{S}f_{1}\wedge\ldots\wedge d_{S}f_{n}=\{f_{1},\ldots,f_{n}\}\omega,
\leqno{(1.7)}$$
where $\omega$ is an $S$-leafwise volume form,
and $d_{S}$ is differentiation along the leaves of $S$. \par
Theorem 1.2 allows us to associate Nambu structures $P$ to the
$n$-dimensional orientable foliations. 
Furthermore,  $\forall f\in
C^{\infty}(M)$, $fP$ also are Nambu structures \cite{DZ} with
singular points at the zeroes of $f$.

It is also essential to stress the following fundamental consequence
of Theorem 1.1 (e.g., \cite{GM}): a Nambu tensor field $P$ necessarily is 
locally decomposable around the regular points and, conversely,
a tensor field $P=V_{1}\wedge ....\wedge V_{n}$ is Nambu iff $span\{V_{1},...,V_{n}\}$
is an involutive distribution on the subset where $P\neq0$.\vspace{2mm}\\
\indent
Given a Nambu bracket of order $n$, if $p$ of its arguments are
fixed, one gets a Nambu bracket of order $n-p$ (a Poisson bracket if
$n-p=2$). The converse is also true \cite{GM}.
On ${\bf R}^{m}$, any constant, decomposable $n$-vector field
$k^{i_{1}\ldots i_{n}}$ is a Nambu tensor \cite{ChT},
and if we keep the fixed function $(1/2)\sum_{j=1}^{m}
(x^{j})^{2}$, we get a Nambu tensor, of order $n-1$,
with the natural components
$$P^{i_{1}\ldots i_{n-1}}=\sum_{j=1}^{m}k^{i_{1}\ldots i_{n-1}j}x^{j}.
\leqno{(1.8)}$$
A Nambu structure defined on a vector space $V^{m}$ by a tensor 
such that its
components with respect to a linear basis of $V$
are linear functions is called a
{\it linear Nambu structure}.
The classification of linear Nambu structures was done by Dufour and Zung
\cite{DZ} (see also \cite{GM}).

Several authors \cite{{Fp},{Tk1},{Gt},{MV},{GM}}
etc. have studied
vector spaces endowed with an
internal, $n$-ary, skew symmetric bracket which satisfies the fundamental
identity of a Nambu bracket.
Following \cite{GM} we call these {\it Filippov algebras}
since they were first studied by Filippov \cite{Fp}.
By looking at brackets of linear functions, it
easily follows  that a linear Nambu structure of order $n$ on
a vector space $V$
induces a Filippov algebra structure on the dual space $V^{*}$.
(The converse may not be true since the structure constants
of a Filippov algebra may form a non decomposable $n$-vector.)

For instance, if $k$ of (1.8) is the canonical volume
tensor of ${\bf R}^{n+1}$ we get the linear Nambu structure of
order $n$ discussed in \cite{ChT}. The corresponding Filippov algebra is the
vector space ${\bf R}^{n+1}$ endowed with the operation of a
{\it vector product}
of $n$ vectors given by
$$v_{1}\times\ldots\times v_{n}=*(v_{1}\wedge\ldots\wedge v_{n}),
\leqno{(1.9)}$$
where $*$ is the Hodge star operator of the canonical Euclidean metric of
${\bf R}^{n+1}$.
The canonical foliation of the linear Nambu structure of
${\bf R}^{n+1}$ defined above has the origin as a $0$-dimensional leaf, and
the spheres with center at the origin as $n$-dimensional leaves. (For
$n=2$,
this is the dual of the Lie algebra $o(3)$ with its well known Lie-Poisson
structure.)\vspace{2mm}

As in the case of Poisson structures, if $(M,P)$ is a
Nambu manifold, and if $p\in M$ is a {\it singular point} of $P$,
the linear part of the Taylor development of $P$ at $p$
defines a linear Nambu structure on $T_{p}M$, and a corresponding
Filippov algebra structure on $T_{p}^{*}M$, which are independent of the
choice of the local coordinates at $p$. This
structure is the
{\it linear approximation} of $P$ at $p$, and $P$ is {\it linearizable} at
$p$ if $P$ is equivalent with its linear approximation on some
neighbourhood of $p$. Linearization theorems were proven in
\cite{DZ}.

This ends the announced recall on Nambu brackets.\vspace{2mm}

Now, following the Poisson model (e.g., \cite{V1}), we will extend
the bracket of functions to a bracket of
$1$-forms. Namely,
for a Nambu structure $P$ of order $n$ on $M^{m}$ we define
$$\{\alpha_{1},\ldots,\alpha_{n}\}=
d(P(\alpha_{(n)}))+\sum_{k=1}^{n}(-1)^{n+k}
i(\sharp_{P}(\alpha_{n,\hat k}))d\alpha_{k}
\leqno{(1.10)}$$
$$=\sum_{k=1}^{n}(-1)^{n+k}L_{\sharp_{P}(\alpha_{n,\hat k})}\alpha_{k}
-(n-1)d(P(\alpha_{(n)})),$$
where $\alpha_{k}$ $(k=1,\ldots,n)$ are $1$-forms on $M$.
The equality of the two expressions of
the new bracket follows from the classical relation
$L_{X}=di(X)+i(X)d$.

The bracket (1.10) will be called the
{\it Nambu form-bracket}, and we have
\proclaim 1.3 Proposition. The Nambu form-bracket satisfies the
following properties: \\i) the form-bracket is totally skew-symmetric;\\
ii) $\forall f_{(n)}\in{\cal F}(M)$, one has
$$\{df_{1},\ldots,df_{n}\}=d\{f_{1},\ldots,f_{n}\};\leqno{(1.11)}$$
iii) for any $1$-forms $\alpha_{(n)}$, and $\forall f\in {\cal F}(M)$ one
has
$$\{f\alpha_{1},\alpha_{2},\ldots,\alpha_{n}\}=
f\{\alpha_{1},\alpha_{2},\ldots,\alpha_{n}\}\leqno{(1.12)}$$
$$+P(df,\alpha_{2},\ldots,\alpha_{n})\alpha_{1}.$$
iv) $\forall f_{(n-1)}\in {\cal F}(M)$ and
for any $1$-form $\alpha$ one has
$$\{df_{1},\ldots,df_{n-1},\alpha\}=L_{X_{f_{(n-1)}}}\alpha.
\leqno{(1.13)}$$
\noindent{\bf Proof.} i) is obvious. ii)
and iii) follow from the first expression
of (1.10).
iv) is a consequence of the first expression (1.10) and of
the commutativity of $d$ and $L$.
Q.e.d.

It would be nice if the form-bracket would also satisfy the fundamental
identity of Nambu brackets. This happens for $n=2$ but, generally,
we only have the following
weaker result
\proclaim 1.4 Proposition. The Lie derivative with respect to a
Hamiltonian vector field is a derivation of the Nambu
form-bracket.\par
\noindent{\bf Proof.} Suppose that the required property holds for the
$1$-forms $\alpha_{(n)}$ i.e.,
$$L_{X_{f_{(n-1)}}}\{\alpha_{1},\ldots,\alpha_{n}\}=
\sum_{k=1}^{n}\{\alpha_{1},\ldots,
L_{X_{f_{(n-1)}}}\alpha_{k},\ldots,\alpha_{n}\}.\leqno{(1.14)}$$
Then, a straightforward computation which uses (1.12) and (1.5) shows that
$L_{X_{f_{(n-1)}}}$ also acts as a derivation of the bracket
$\{f\alpha_{1},\alpha_{2},\ldots,\alpha_{n}\}$,
$\forall f\in{\cal F}(M)$.

This remark shows that the proposition is true if the
result holds for a bracket of the form
$\{dg_{1},\ldots,dg_{n}\}$, $\forall g_{k}\in{\cal F}(M)$.
But, this  is a consequence of the fundamental identity for functions
since by (1.11) we have
$$L_{X_{f_{(n-1)}}}\{dg_{1},\ldots,dg_{n}\}=
L_{X_{f_{(n-1)}}}d\{g_{1},\ldots,g_{n}\}=
dL_{X_{f_{(n-1)}}}\{g_{1},\ldots,g_{n}\}.$$
Q.e.d.

The relation between (1.14) and the fundamental identity for $1$-forms is
given by (1.13). Since any closed form locally is an exact form, we
see that the {\it fundamental identity}
$$\{\beta_{1},\ldots,\beta_{n-1},
\{\alpha_{1},\ldots,\alpha_{n}\}\}=
\sum_{k=1}^{n}\{\alpha_{1},\ldots,\alpha_{k-1},\leqno{(1.15)}$$
$$\{\beta_{1},\ldots,\beta_{n-1},\alpha_{k}\},\alpha_{k+1},\ldots,
\alpha_{n}\}$$
{\it holds whenever the} $1$-{\it forms} $\beta$ {\it are closed}.

In \cite{GM}, the bracket (1.10) was interpreted
as a dual bracket of the complete lift of $P$ to $TM$ .
\section{Nambu-Lie Groups}
Since Poisson-Lie groups play an important role in Poisson geometry (e.g.,
\cite{V1}), we are motivated to discuss similarly defined Nambu-Lie
groups. These cannot be defined by the demand that the multiplication is a
Nambu morphism since, on a product manifold, the sum of Nambu tensors
may not be a Nambu tensor. But, it makes sense to say that
a Nambu tensor $P$ of order $n$ endows the Lie group $G$ with the structure of a
{\it Nambu-Lie group} if $P$ is a {\it multiplicative tensor field}
i.e. (e.g., \cite{V1}), $\forall g_{1},g_{2}\in G$, one has
$$P_{g_{1}g_{2}}=L_{g_{1}^{*}}P_{g_{2}}+R_{g_{2}^{*}}P_{g_{1}},
\leqno{(2.1)}$$
where $L$ and $R$ denote left and right translations in $G$, respectively.

The multiplicativity of $P$ implies $P(e)=0$, where $e$ is the unit of $G$.
Moreover, if $G$ is connected, $P$ is multiplicative iff $P(e)=0$, and the
Lie derivative $L_{X}P$ is a left (right) invariant tensor field whenever
$X$ is left (right) invariant. As an immediate
consequence it follows that
the Nambu-Lie group structures on the additive Lie group
${\bf R}^{m}$ are exactly the linear Nambu structures
of ${\bf R}^{m}$.

From (2.1), it follows that the set
$$G_{0}:=\{g\in G\;/\;P_{g}=0\}$$
is a closed subgroup of $G$. Indeed,
(2.1) shows that if $g_{1},g_{2}\in
G_{0}$, the product $g_{1}g_{2}\in G_{0}$. Furthermore,
if $g\in G_{0}$, then
$$0=P_{e}=P_{gg^{-1}}=L_{g^{*}}P_{g^{-1}},$$
hence $g^{-1}\in G_{0}$. 

The following theorem extends the characterization of the Poisson-Lie
groups given by
Dazord and Sondaz \cite{DzS}.
\proclaim 2.1 Theorem. If $G$ is a connected Lie group endowed with a
Nambu tensor field $P$ which vanishes at the unit $e$ of $G$, then
$(G,P)$ is a Nambu-Lie group iff the $P$-bracket of any $n$ left
(right) invariant $1$-forms of $G$ is a left (right) invariant $1$-form.
\par
\noindent{\bf Proof.} As in the Poisson case (e.g.,
\cite{V1}), the evaluation of the Lie derivative
via (1.10) yields
$$(L_{Y}\{\alpha_{1},\ldots,\alpha_{n}\})(X)=
Y((L_{X}P)(\alpha_{(n)}))$$
for any left invariant vector field $X$, right invariant vector field $Y$,
and left invariant $1$-forms $\alpha_{(n)}$. (Same if left and right are
interchanged.) Hence, the condition of the theorem is equivalent with
the fact that $L_{X}P$ is left invariant whenever $X$ is left invariant. Q.e.d.
\vspace{2mm}

Many other properties of Poisson-Lie groups also have a
straightforward generalization.

Since $P(e)=0$,
the linear approximation of $P$ at $e$ defines a linear Poisson structure
on the Lie algebra ${\cal G}$ of $G$ and a dual Filippov algebra structure
on the dual space ${\cal G}^{*}$. Furthermore,
as for $n=2$, a compatibility relation between
the Lie bracket and the linear Nambu structure of ${\cal G}$
exists.

Following \cite{LW},
let us consider the {\it intrinsic derivative} $\pi_{e}:=d_{e}P:
{\cal G}\rightarrow\wedge^{n}{\cal G}$ defined by
$$\pi_{e}(X)(\alpha_{(n)})=(L_{\check X}P)_{e}(\alpha_{(n)})
,\leqno{(2.2)}$$
where $\alpha_{(n)}\in{\cal G}^{*}$, $X\in{\cal G}$, and $\check X$ is any
vector field on {\cal G} with the value $X$ at $e$. 
It is easy to understand that, if $\Pi$ is the Nambu tensor on ${\cal G}$
associated with the linear approximation of $P$ at $e$, then, 
$\forall X\in{\cal G}$, $\Pi_{X}=\pi_{e}(X)$.
Furthermore, we have
\proclaim 2.2 Theorem. i). The bracket of the dual Filippov algebra
structure of ${\cal G}^{*}$ is the dual
$\pi^{*}_{e}$ of the mapping $\pi_{e}$, and it
has each of the following expressions
$$[\alpha_{1},\ldots,\alpha_{n}]=
d_{e}(P(\check\alpha_{(n)}))=\pi_{e}^{*}(\wedge^{n}_{k=1}\alpha_{k})
\leqno{(2.3)}$$
$$=\{\bar\alpha_{1},\ldots,\bar\alpha_{n}\}_{e}
=\{\tilde\alpha_{1},\ldots,\tilde\alpha_{n}\}_{e},$$
where $\alpha_{(n)}\in{\cal G}^{*}$, $\check\alpha_{(n)}$ are $1$-forms on $G$
which are equal to $\alpha_{(n)}$ at $e$, and $\bar\alpha_{(n)}$,
$\tilde\alpha_{(n)}$ are the left and right invariant $1$-forms,
defined by $\alpha_{(n)}$, respectively.\\
ii). The mapping $\pi_{e}$ is a $\wedge^{n}{\cal
G}$-valued $1$-cocycle of ${\cal G}$
with respect to the adjoint representation
$$ad_{X}(Y_{1}\wedge\ldots\wedge Y_{n})=
\sum_{k=1}^{n}Y_{1}\wedge\ldots Y_{k-1}
\wedge[X,Y_{k}]_{{\cal G}}\wedge Y_{k+1}
\wedge\ldots\wedge Y_{n},$$
($X,Y_{(n)}\in{\cal G}$). \par
\noindent{\bf Proof.} The proofs are as for $n=2$; see
\cite{LW} or Chapter 10 of
\cite{V1}.

i). The first equality sign is just the definition
of the linear approximation,
and the second follows since from $P_{e}=0$ we have
$$<\pi_{e}^{*}(\wedge^{n}_{k=1}\alpha_{(k)}),X>=\pi_{e}(X)(\alpha_{(n)})
=(L_{\check X}P)_{e}(\alpha_{(n)})=X(P(\check
\alpha_{(n)}))$$ $$=<d_{e}(P(\check\alpha_{(n)}),X>.$$
The remaining part of (2.3) is proven by the following computation with
left (similarly, right) invariant forms:
$$\{\bar\alpha_{1},\ldots,\bar\alpha_{n}\}_{e}(X)\stackrel{(1.10)}{=}
X(P(\bar\alpha_{(n)}))+\sum_{k=1}^{n}
(-1)^{n+k}(d\bar\alpha_{k})_{e}(\sharp_{P}(\alpha_{(n,\hat
k)}),X)$$
$$=X(P(\bar\alpha_{(n)}))-\sum_{k=1}^{n}(-1)^{n+k}
(L_{\tilde X}\bar\alpha_{k})_{e}(\sharp_{P}(\alpha_{(n,\hat
k)}))$$
$$+\sum_{k=1}^{n}\sharp_{P}(\alpha_{(n,\hat
k)})_{e}(\bar\alpha_{k}(\tilde X))=X(P(\bar\alpha_{(n)}),$$
where $\tilde X$ is the right invariant vector field defined by $X$,
and we used the equalities $P_{e}=0$, $L_{\tilde X}\bar\alpha_{k}=0$.

ii). The fact that $\pi_{e}$ is a $1$-cocycle means that we have
$$\partial\pi_{e}(X,Y):=
ad_{X}(\pi_{e}(Y))-ad_{Y}(\pi_{e}(X))-\pi_{e}([X,Y]_{{\cal G}})=0
,\leqno{(2.4)}$$
where $X,Y\in {\cal G}$,
or equivalently,
$$ad_{X}(\Pi_{Y})-ad_{Y}(\Pi_{X})=\Pi_{[X,Y]}, \leqno{(2.4')}$$
with the earlier defined $\Pi$.
We always use the notation with bars and tildes for left and right invariant
objects on Lie groups as we did above. Then,
$$ad_{X}(\pi_{e}(Y))=\frac{d}{ds}/_{s=0}Ad\;exp(sX)((L_{\bar Y}P)_{e})=
(L_{\bar X}L_{\bar Y}P)_{e},$$
and (2.4) is a consequence of this result. Q.e.d.

The second part of Theorem 2.2 yields
\proclaim 2.3 Corollary. $\forall\alpha_{(n)}\in{\cal G}^{*}$
and $\forall X,Y\in{\cal G}$ the following
relation holds
$$<[\alpha_{1},\ldots,\alpha_{n}],[X,Y]_{{\cal G}}>=
\sum_{k=1}^{n}(<[\alpha_{1},\ldots,\alpha_{k-1},coad_{Y}\alpha_{k},
\leqno{(2.5)}$$
$$\alpha_{k+1},\ldots,\alpha_{n}],X>
-<[\alpha_{1},\ldots,\alpha_{k-1},coad_{X}\alpha_{k},
\alpha_{k+1},\ldots,\alpha_{n}],Y>).$$
\par
\noindent{\bf Proof.} The result is nothing but a reformulation of the
cocycle condition (2.4). Q.e.d. 

On the other hand, the first part of
Theorem 2.2  allows us to get a result on subgroups just as in the 
Poisson case. A Lie subgroup $H$ of a Nambu-Lie group $(G,P)$ will be called 
a {\it Nambu-Lie subgroup} if $H$ has a (necessarily unique) multiplicative 
Nambu tensor $Q$ such that $(H,Q)$ is a Nambu submanifold of $(G,P)$.
For instance, the vanishing subgroup $G_{0}$ of $P$ 
with the Nambu structure $Q=0$
is a Nambu-Lie subgroup
of $(G,P)$.
If $H$ is connected, it is a Nambu-Lie subgroup of $(G,P)$
iff $Ann({\cal H})$, where ${\cal H}$ is the Lie algebra of $H$,
is an ideal in $({\cal G}^{*},[.,...,.])$. By this we mean that the bracket
(2.3) is in $Ann({\cal H})$ whenever one of the arguments (at least)
is in $Ann({\cal H})$. The proof is the same as for $n=2$ e.g., \cite{V1}.

Furthermore, if $(H,Q)$ is a Nambu-Lie subgroup of $(G,P)$,
the homogeneous space $M=G/H$ inherits a Nambu structure $S$ of the 
same order as $P,Q$ such that the natural projection
$p:(G,P)\rightarrow(M,S)$ is a Nambu morphism. This holds since the 
brackets  $\{f_{1}\circ p,...,f_{n}\circ p\}_{P}$ 
$(f_{i}\in C^{\infty}(M))$ are
constant along the fibers of $p$, which is easy to check using (2.1).
(E.g., see Proposition 10.30 in \cite{V1} for the case $n=2$.)
Moreover, as a consequence of (2.1), the natural left action
of $G$ on $M$ satisfies the multiplicativity condition
$$S_{g(x)}=\varphi_{g^{*}}(S_{x})+\varphi^{x}_{*}(P_{g}), \leqno{(2.1')}$$
where $g\in G,x\in M$, and
$\varphi_{g}:M\rightarrow M,\: \varphi^{x}:G\rightarrow M$ are defined by
$\varphi_{g}(x)=\varphi^{x}(g)=g(x)$.
An interesting example is $H=G_{0}$, $Q=0$, $M=G/G_{0}$,
where $G_{0}$ is the vanishing subgroup of $P$.
In agreement with the above situation, any action of a 
Nambu-Lie group $(G,P)$ on a Nambu manifold $(M,S)$ which
satisfies (2.1$'$) will be 
called a {\it Nambu action}. If $G$ is connected, one has the
same infinitesimal characteristic properties of Nambu actions 
as in the Poisson case e.g., Proposition 10.27 in \cite{V1}.
In particular, that $\forall X\in{\cal G}$,
$L_{X_{M}}S=-[(d_{e}P)(X)]_{M}$, where $e$ is the unit of $G$, and 
the index $M$ denotes the infinitesimal action on $M$.
\vspace{2mm}\\ 
\indent
At this point, one may ask whether a Lie algebra ${\cal G}$
with a linear Nambu structure $\Pi$ that satisfies the cocycle condition (2.4), (2.5)
can be integrated to a Nambu-Lie group $(G,P)$?

Some of the results
known for $n=2$ still hold.
If $G$ is connected and simply
connected,
for any $1$-cocycle $\pi_{e}$ as in Theorem 2.2 ii),
there exists a unique multiplicative $n$-vector field $P$ on $G$,
called the {\it integral field} of $\pi_{e}$ such that
$d_{e}P$ is the given cocycle.
Indeed, for the given cocycle $\pi_{e}$, $$\pi_{g}(X_{g}):=Ad\,g(
\pi_{e}(L_{g^{-1*}}X_{g}))\hspace{5mm}(g\in G,\;X_{g}\in T_{g}G)$$
defines a $\wedge^{n}{\cal G}$-valued $1$-form
$\pi$ on $G$ which satisfies the
equivariance condition $L_{g}^{*}\pi=(Ad\,g)\circ\pi$. This implies that
$d\pi=0$, and, since $G$ is connected and simply connected, $\pi=dP$ for a
unique $n$-vector field $P$ on $G$, which can be seen to be
multiplicative \cite{{LW},{V1}}.
If this field is Nambu, we are done. But, this final part
is not ensured if $n\geq3$.

Moreover, the structure theory of multiplicative Nambu tensors of \cite{GM}
leads to one more, important, necessary condition.
Recalling from \cite{GM}, let $(G,P)$ be a connected Nambu-Lie group,
and $\forall g\in G$ where $P_{g}\neq0$ put $\tilde P(g)=R_{g^{-1*}}P_{g}$.
$\tilde P(g)$ is a decomposable element of $\wedge^{n}{\cal G}$, and
yields a subspace $V(g):=$ {\mit span of the factors of} $\tilde P(g)\subseteq{\cal G}$.
Furthermore, $V_{\cap}:=\cap_{g}V(g)$, and $V_{\cup}:=\sum_{g}V(g)$
are ideals in ${\cal G}$, 
and a sum-intersection lemma (\cite{DZ}, Lemma 3.2 and Theorem 3.1,
\cite{GM}, Lemma 1) tells that either $dim\,V_{\cap}\geq n-1$ or
$dim\,V_{\cup}= n+1$.
In the first case, put ${\cal H}=V_{\cap}$, and if $dim\, V_{\cup}=n+1$
put ${\cal H}=V_{\cup}$. In both cases, the ideal ${\cal H}$
will be called the {\it core ideal} of $P$, and if $0\neq\Lambda_{0}\in
\wedge^{dim\,{\cal H}}{\cal H}$, we call $\Lambda_{0}$ a {\it core} of $P$.
According to the cases $dim\,{\cal H}=n,n-1,n+1$ we have 
either a) $\tilde P(g)=\theta(g)\Lambda_{0}$,
with $\theta\in C^{\infty}(G)$, or b) $\tilde P(g)=X(g)\wedge\Lambda_{0}$
with $X:G\rightarrow{\cal G}$, or c) $\tilde P(g)=i(\alpha(g))\Lambda_{0}$
with $\alpha:G\rightarrow{\cal G}^{*}$, and the conditions which $\theta,X,\alpha$ 
must satisfy in order to provide a multiplicative Nambu tensor $P$ on $G$
are determined in \cite{GM}. In case \hspace{1mm}
a) the canonical foliation ${\cal D}$ of $P$ is
the same as the left (and right) invariant foliation 
${\cal F}_{\cal H}$ defined by translating ${\cal
H}$ along $G$, in case b) ${\cal F}_{\cal H}$ is a subfoliation of codimension $1$ of
${\cal D}$, and in case c) ${\cal D}$ is a subfoliation of codimension $1$ of ${\cal
F}_{\cal H}$.

In particular, since any linear Nambu structure $\Pi$ on a vector space $V$
is multiplicative with respect to the additive structure of $V$,
$\Pi$ has a core linear subspace ${\cal H}$ (ideal of a commutative Lie algebra) 
and a core $\Lambda_{0}$. Moreover, we get
\proclaim 2.4 Theorem. Let $(G,P)$ be a connected,
$m$-dimensional, Nambu-Lie group with the 
core ideal ${\cal H}$ and a core $\Lambda_{0}$,
and let $\Pi$ be the linear approximation of $P$ at the unit $e\in G$.
Then, ${\cal H},\Lambda_{0}$ also are the core subspace and a core of $\Pi$, 
respectively. Furthermore, the core of $\Pi$ is an ideal of ${\cal G}$,
and ${\cal G}$ must have an ideal ${\cal H}$ of one of the dimensions $n,n-1,n+1$,
where $n$ is the order of $P$.
\par
\noindent{\bf Proof.}
By its very definition, $\pi_{X}=\pi_{e}(X)=X\tilde P$, where $\tilde P$ is seen as a
$\wedge^{n}{\cal G}$-valued function on $G$ (e.g., \cite{V1}, p. 166). Then, by
derivating the three possible expressions of $\tilde P(g)$ as recalled from \cite{GM}
above, we get the conclusion. Q.e.d.

Theorem 2.4 reduces the finding of the core ideal and the core of a Nambu-Lie
group to the same problem for a linear Nambu structure.

Now, in agreement with
Theorems 2.2 and 2.4, we define a {\it Nambu-Lie
algebra} as being a Lie algebra ${\cal G}$ endowed with a linear Nambu structure 
$\Pi$ such that: i) $\Pi$ is a $1$-cocycle , and ii) the core subspace of $\Pi$
is an ideal of ${\cal G}$.

The method of \cite{GM} can 
also be used for the determination of the Nambu-Lie algebras $({\cal G},\Pi)$.
Namely, if ${\cal G}$ is given, we have to consider the ideals ${\cal H}$
of ${\cal G}$. Then, for an ideal ${\cal H}$ of dimension $n$, and with
$0\neq\Lambda_{0}\in\wedge^{n}{\cal H}$, we have to look for $\Pi_{X}$ under one of
the 
following forms:
$$\Pi_{X}=\varphi(X)\Lambda_{0},\hspace{1cm}\varphi\in{\cal G}^{*},\leqno{(2.6)}$$
$$\Pi_{X}={\cal X}(X)\wedge\Lambda_{0},\hspace{1cm}{\cal X}\in End({\cal G}),\leqno{(2.7)}$$
$$\Pi_{X}=i(\alpha(X))\Lambda_{0},
\hspace{1cm}\alpha\in Hom({\cal G},{\cal G}^{*}).\leqno{(2.8)}$$
Finally, we must ask the cocycle condition (2.4$'$) which,
respectively, will give
$$\varphi([X,Y])=\varphi(Y)\gamma(X)-\varphi(X)\gamma(Y),\leqno{(2.6')}$$
$${\cal X}([X,Y])-[X,{\cal X}(Y)]+[Y,{\cal X}(X)]-
\gamma(X){\cal X}(Y)+\gamma(Y){\cal X}(X)\in{\cal H}, \leqno{(2.7')}$$
$$i[\alpha[X,Y])+(coad_{X}\alpha)(Y)-(coad_{Y}\alpha)(X) \leqno{(2.8')}$$
$$+\gamma(X)\alpha(Y)-\gamma(Y)\alpha(X)]\Lambda_{0}=0,$$
where $X,Y\in{\cal G}$, $coad$ is naturally extended to the twice covariant tensor
$\alpha$,
and $\gamma\in {\cal G}^{*}$
is determined by the condition
$ad_{X}\Lambda_{0}=\gamma(X)\Lambda_{0}$, therefore, it must satisfy
$$\gamma([X,Y])=0,\hspace{1cm}\forall X,Y\in {\cal G}. \leqno{(2.9)}$$

For an example, let us consider the unitary Lie algebra $u(2)$.
In $u(2)$ we have the basis
$$X_{1}=\frac{\sqrt{-1}}{2}\left(\begin{array}{rr}1&0\\0&1\end{array}\right),\;
X_{2}=-\frac{\sqrt{-1}}{2}\left(\begin{array}{rr}0&1\\1&0\end{array}\right),\;$$
$$X_{3}=-\frac{1}{2}\left(\begin{array}{rr}0&1\\-1&0\end{array}\right),\;
X_{4}=-\frac{\sqrt{-1}}{2}\left(\begin{array}{rr}1&0\\0&-1\end{array}\right),$$
such that $X_{1}$ spans the center, and
$$[X_{2},X_{3}]=X_{4},\;[X_{3},X_{4}]=X_{2},\;[X_{4},X_{2}]=X_{3},\leqno{(2.10)}$$
and we will denote by $(x^{a})$ $(a=1,2,3,4)$, 
the corresponding linear coordinates.The relevant ideals 
are ${\cal H}_{1}=u(2)$, ${\cal H}_{2}=su(2)=span\{X_{2},X_{3},X_{4}\}$,
with the cores
$$\Lambda_{01}=X_{1}\wedge X_{2}\wedge X_{3}\wedge X_{4},\;
\Lambda_{02}=X_{2}\wedge X_{3}\wedge X_{4},$$ respectively.

From (2.9), (2.10), we see that $\gamma=0$ ,
and any $\varphi$ which satisfies (2.6$'$) vanishes at $X_{2},X_{3},X_{4}$.
It follows that on $u(2)$ there is only one Nambu-Lie
algebra structure of type (2.6) with the core $\Lambda_{01}$, up to a constant factor, 
and this is
$$\Pi=x^{1}\frac{\partial}{\partial x^{1}}
\wedge\frac{\partial}{\partial x^{2}}
\wedge\frac{\partial}{\partial x^{3}}
\wedge\frac{\partial}{\partial x^{4}}.\leqno{(2.11)}$$

But, it is easy to see that there is no non zero Nambu-Lie group structure
$P$ on the unitary group
$U(2)$ with the core $\Lambda_{01}$. Such a structure would have
$\tilde P(g)=\theta(g)\Lambda_{01}$ where multiplicativity implies
that $\forall g_{1},g_{2}\in U(2)$ one has
$\theta(g_{1}g_{2})=\theta(g_{1})+\theta(g_{2})$ \cite{GM}. 
And, there is no non zero $\theta$ with this property
since by
Theorem 6 of \cite{GM} one should have $\theta(unit)=0$, $d\theta=$ a bi-invariant
$1$-form, i.e., $d\theta=kdx^{1}$ ($k=const.$). Hence, $\theta$ would be an additive
character
of the circle subgroup $S^{1}$ of $U(2)$ and, thus,  $\theta=0$.

Therefore, above we have an example of a Nambu-Lie algebra which does not integrate
to a Nambu-Lie group.

Furthermore, for $\Lambda_{01}$, $\Pi$ of (2.7) is zero. But, it is 
possible to find Nambu-Lie structures of the type (2.8). An example is
$$\Pi=(x^{2}\frac{\partial}{\partial x^{2}}+x^{4}
\frac{\partial}{\partial x^{4}})\wedge\frac{\partial}{\partial x^{3}}
\wedge\frac{\partial}{\partial x^{1}}=\partial(X_{4}\wedge X_{2}\wedge X_{1}),
\leqno{(2.12)}$$
which is a coboundary hence, a cocycle.
Here, $\partial X_{1}=0$ since $X_{1}$ is $ad$-invariant, and $\partial
(X_{4}\wedge X_{2})$ is the cocycle of a well known example
of a Poisson-Lie structure of $SU(2)$ namely \cite{Lu},
$$W(g)=L_{g^{*}}(X_{4}\wedge X_{2})-R_{g^{*}}(X_{4}\wedge X_{2}),\;\;(g\in SU(2)).
\leqno{(2.13)}$$
Accordingly, (2.12) is the cocycle of a Nambu-Lie structure on $U(2)$
which is defined by
$$P(g)=W(g)\wedge X_{1} \leqno{(2.14)}$$
$$=L_{g^{*}}(X_{4}\wedge X_{2}\wedge X_{1})-
R_{g^{*}}(X_{4}\wedge X_{2}\wedge X_{1}),\;\;(g\in U(2)).$$
Indeed, it is easy to check that $P$ is multiplicative \cite{V1}.
It is decomposable since,
if $W(g)\neq0$, $rank\,W(g)=2$, and the factors span an involutive distribution.

It is interesting that we have thereby obtained an example of a compact
Nambu-Lie group. The construction does not extend to $U(n)$ with $n>2$ since the 
similar structure $P$ is not decomposable.

Now, comming back to the Lie algebra $u(2)$, we should look at
cocycles with the core ideal ${\cal H}_{2}=su(2)$. Like for (2.11),
it follows again that, up to a constant factor, the only structure of type (2.6) is
$$P=x_{1}\frac{\partial}{\partial x_{2}}\wedge
\frac{\partial}{\partial x_{3}}\wedge\frac{\partial}{\partial
x_{4}},\leqno{(2.15)}$$
and it does not come from a Nambu-Lie group structure on $U(2)$. The structures (2.7) 
which can be obtained reduce to (2.11), and (2.8) may 
lead to Poisson-Lie structures on $U(2)$.
In particular, we can refind (2.13) with $g\in U(2)$.
\vspace{1mm}\\ \indent
In principle, a systematic search for the structures $\Pi$ should be possible by 
using the canonical forms of the linear Nambu structures given in \cite{{DZ},{GM}}.
One has to look for structure constants of Lie algebras
which, together with the canonical structures of \cite{DZ}, satisfy the 
cocycle condition (2.5), and such that the core of the linear Nambu structure is an
ideal of the considered Lie algebra.\vspace{1mm}\\
\indent
We finish by giving some more
examples of  non commutative Nambu-Lie
groups. 

A first example is the $3$-dimensional solvable Lie group
$$G_{3}:=\{\left( \begin{array}{rrr} x&0&y\\ 0&x&z\\0&0&1 \end{array}
\right)\;\;/x,y,z\in {\bf R},\;x\neq0\}.\leqno{(2.16)}$$
The left invariant forms of this group are $dx/x,dy/x,dz/x$, and if we look
for a Nambu tensor of the form
$$P=f(x)\frac{\partial}{\partial x}\wedge\frac{\partial}{\partial y}\wedge
\frac{\partial}{\partial z}\leqno{(2.17)}$$
such that $\{dx/x,dy/x,dz/x\}$ is left-invariant, and $f(1)=0$, we see that
$f=x(x^{2}-1)/2$ does the job. The corresponding Nambu-Lie algebra is ${\bf
R}^{3}$ with the linear Nambu structure $x^{1}(\partial/\partial x^{1})
\wedge(\partial/\partial x^{2})\wedge(\partial/\partial x^{3})$.

The next example is the {\it generalized Heisenberg group}
$$H(1,p):=\{\left(\begin{array}{rrr}
Id_{p}&X&Z\\ 0&1&y\\0&0&1 \end{array}\right)\},\leqno{(2.18)}$$
where $X=\,^{t}(x_{1}...x_{p})$, $Z=\,^{t}(z_{1}...z_{p})$
($t$ means transposition of matrices).
The left invariant $1$-forms of this group are
$$dx_{1},...,dx_{p},dy,dz_{1}-x_{1}dy,...,dz_{p}-x_{p}dy,\leqno{(2.19)}$$
and $$P=y\frac{\partial}{\partial x_{1}}\wedge
\frac{\partial}{\partial z_{1}}\wedge\frac{\partial}{\partial y}
\leqno{(2.20)}$$
makes $H(1,p)$ into a Nambu-Lie group. Indeed, it vanishes at the unit, and
it follows easily that the brackets of the left invariant $1$-forms are
left invariant. The corresponding Nambu-Lie algebra is ${\bf R}^{2p+1}$
with the same Nambu tensor (2.20).

A third example is the direct product
$G=H(1,1)\times{\bf R}_{+}$, where ${\bf R}_{+}$ is the multiplicative
group of the positive real numbers $t$. The left invariant $1$-forms
of the group are those given by (2.19), and $dt/t$. The tensor
$$P=t(\ln t)\frac{\partial}{\partial y}\wedge
\frac{\partial}{\partial z}\wedge\frac{\partial}{\partial
t}\leqno{(2.21)}$$
makes $G$ into a Nambu-Lie group for the same reasons as in the previous
examples. The corresponding Nambu-Lie algebra is ${\bf R}^{4}$ with the
linear Nambu structure
$P=x_{4}(\partial/\partial x_{2})\wedge
(\partial/\partial x_{3})\wedge(\partial/\partial
x_{4})$.

We also notice that if $(G_{1},P)$ is a Nambu-Lie group, and $G_{2}$ is
any other Lie group, $fP$, where $f\in C^{\infty}(G_{2})$, is a Nambu-Lie
structure on $G_{1}\times G_{2}$ (check (2.1)).

Finally, we quote the following 
important result proven in \cite{GM}: there
are no Nambu-Lie structures of order $n\geq3$ on simple Lie groups,
and if $G=G_{1}\times...\times G_{s}$ is a semisimple Lie group
with the simple factors $G_{i}$ $(i=1,...,s)$, the only multiplicative
Nambu tensors on $G$ are wedge products of contravariant volume tensor fields
on a part of the factors with either multiplicative Poisson bivectors
or multiplicative vector fields on other factors.
 \vspace{5mm}
{\small Department of Mathematics, \\}
{\small University of Haifa, Israel.\\}
{\small E-mail: vaisman@math.haifa.ac.il\\}
\end{document}